\documentclass[11pt,leqno,oneside]{article}
\usepackage{amsmath}
\usepackage{amsthm}
\usepackage{amstext}
\usepackage{amsopn}
\usepackage{amsbsy}

\newtheorem*{main-theorem}{Main Theorem}
\newtheorem{theorem}{Theorem}

\newtheorem{lemma}{Lemma}

\begin{document}

\title{A Lower Bound of the First Dirichlet Eigenvalue of a Compact Manifold
with Positive Ricci Curvature
\thanks{2000 Mathematics Subject Classification Primary 58J50, 35P15; Secondary 53C21}}

\author{Jun LING}
\date{}

\maketitle

\begin{abstract}
We give a new estimate on the lower bound for the first Dirichlet
eigenvalue for a compact manifold with positive Ricci curvature in
terms of the in-diameter and the lower bound of the Ricci
curvature. The result improves the previous estimates.
\end{abstract}

%%%%%%%%%%%%%%%%%%Introduction%%%%%%%%%%%%%%%%%%
\section{Introduction}\label{sec-intro}
If $(M,g)$ is an n-dimensional compact Riemannian manifold whose
Ricci curvature has a positive lower bound $(n-1)K$ for some
constant $K>0$ and whose non-empty boundary $\partial M$ has
nonnegative mean curvature with respect to the outward normal,
Reilly \cite{reilly} gave the following lower bound of the first
Dirichlet eigenvalue $\lambda$ of the Laplacian on $M$
\begin{equation}\label{nK-bound}
\lambda\geq nK.
\end{equation}
This estimate gives no information when the above constant $K$
vanishes. In such case, Li-Yau \cite{liy1} and Zhong-Yang
\cite{zy} provided another lower bound for the first non-zero
eigenvalue of a closed manifold
\[
\lambda\geq \frac{\pi^2}{d^2}.
\]

It is an interesting problem to find a unified lower bound of the
first Dirichlet eigenvalue $\lambda$ in terms of the lower bound
$(n-1)K$ of the Ricci curvature and the diameter $d$, in-diameter
$\tilde{d}$ and other geometric quantities, which do not vanish as
$K$ vanishes, of the manifold with positive Ricci curvature. D.
Yang \cite{yang} proved that
\begin{equation}\label{yang-bound}
\lambda \geq \frac14(n-1)K + \frac{\pi^2}{(\tilde{d})^2},
\end{equation}
where $\tilde{d}$ is the diameter of the largest interior ball in
$M$.

In this paper we give a new estimate on the lower bound of the
first Dirichlet eigenvalue $\lambda$. We have the following
result.

\begin{theorem} \label{main-thm}
If $(M,g)$ is an $n$-dimensional compact Riemannian manifold with
boundary. Suppose that Ricci curvature \textup{Ric}$(M)$ of $M$ is
bounded below by $(n-1)K$ for some constant $K>0$
\begin{equation}\label{ricci-bound}
\textup{Ric}(M)\geq (n-1)K
\end{equation}
and that the mean curvature of the boundary $\partial M$ with
respect to the outward normal is nonnegative, then the first
Dirichlet eigenvalue $\lambda$ of the Laplacian $\Delta$ of $M$
has the following lower bound
\begin{equation}                    \label{main-bound-in}
\lambda \geq  \frac12(n-1)K+\frac{\pi^2}{(\tilde{d})^2},
\end{equation}
where $\tilde{d}$ is the diameter of the largest interior ball in
$M$, that is, $\tilde{d}=2\sup_{x\in M}\{\textup{dist}(x,\partial
M)\}$.
\end{theorem}

Our result improves Yang' bound (\ref{yang-bound}) by doubling the
coefficient before $(n-1)K$. In the proof,  we use a function
$\xi$ that the author constructed in \cite{ling-d} for the
construction of the suitable test function instead of using the
Zhong-Yang's canonical function. That provides a new way to
sharpen the bound. In the next section, we derive some preliminary
estimates and conditions for test functions  first and we
construct the needed test function and prove the main result in
the last section.

%%%%%%%%%%%%%Estimates and Barriers%%%%%%%%%%%%%%%
\section{Preliminary Estimates}\label{sec-pre-es}
The first basic estimate is of Lichnerowicz-type. Recall that the
classic Lichnerowicz Theorem \cite{lich} states that if $M$ is an
$n$-dimensional closed manifold whose Ricci curvature satisfies
(\ref{ricci-bound}) then the first non-zero eigenvalue has a lower
bound (\ref{nK-bound}). Reilly \cite{reilly} proved that this
Lichnerowicz-type estimate remains true for the first Dirichlet
eigenvalue $\lambda$ as well if the manifold has the same lower
bound for the Ricci curvature and has non-empty boundary whose
mean curvature with respect to the outward normal is nonnegative.
For the completeness and consistency, we use gradient estimate in
\cite{li}-\cite{liy1} and \cite{sy} to derive the
Lichnerowicz-type estimate.
\begin{lemma}               \label{nK-lemma}
Under the conditions in Theorem \ref{main-thm}, the estimate
(\ref{nK-bound}) holds.
\end{lemma}
\begin{proof}\quad
Let $v$ be a normalized eigenfunction of the first Dirichlet
eigenvalue such that
\begin{equation}                                \label{v-d-con}
\sup_{M}v=1, \quad \inf_{M}v=0.
\end{equation}
The function $v$ satisfies the following
\begin{equation}                                 \label{v-d-eq}
\Delta v=-\lambda v\quad \textrm{in }M
\end{equation}
and
\begin{equation}                                \label{v-d-boundary}
v=0\qquad \textrm{ on }\partial M.
\end{equation}

Take an orthonormal frame $\{e_1, \dots, e_n\}$ of $M$ about
$x_0\in M$. At $x_0$ we have
\[
\nabla_{e_j}(|\nabla v|^2)(x_0)= \sum_{i=1}^n2v_iv_{ij}
\]
and
\begin{eqnarray}
\Delta (|\nabla v|^2)(x_0)&=&2\sum_{i,j=1}^nv_{ij}v_{ij}+2\sum_{i,j=1}^nv_{i}v_{ijj}\nonumber\\
      &=&2\sum_{i,j=1}^nv_{ij}v_{ij}+2\sum_{i,j=1}^nv_{i}v_{jji}+2\sum_{i,j=1}^n\textrm{R}_{ij}v_iv_j\nonumber\\
      &=&2\sum_{i,j=1}^nv_{ij}v_{ij}+2\nabla v \nabla (\Delta v) + 2\textrm{Ric}(\nabla v, \nabla v)\nonumber\\
      &\geq&2\sum_{i=1}^nv_{ii}^2+ 2\nabla v \nabla (\Delta v) + 2(n-1)K|\nabla
      v|^2\nonumber\\
      &\geq&\frac{2}{n}(\Delta v)^2-2\lambda|\nabla v |^2 + 2(n-1)K|\nabla v|^2.\nonumber
\end{eqnarray}
Thus at all point $x\in M$,
\begin{equation}                        \label{basic2.002-d}
\frac12\Delta (|\nabla v|^2)\geq \frac1n \lambda^2 v^2 +
[(n-1)K-\lambda]|\nabla v|^2.
\end{equation}
On the other hand, after multiplying (\ref{v-d-eq}) by $v$ and
integrating both sides over $M$ and using (\ref{v-d-boundary}), we
have
\[
\int_{M}\lambda v^2\,dx=-\int_Mv\Delta v\, dx
\]
\[=-\int_{\partial M}v\frac{\partial}{\partial
\nu}v\,ds+\int_{M}|\nabla v|^2\, dx=\int_{M}|\nabla v|^2\, dx,
\]
where and below $\nu$ is the outward normal of $\partial M$. That
the integral on the boundary vanishes is due to
(\ref{v-d-boundary}). Integrating (\ref{basic2.002-d}) over $M$
and using the above equality, we get
\begin{equation}                                        \label{basic2.003-d}
\frac12\int_{\partial M}\frac{\partial}{\partial \nu}(|\nabla
v|^2)\,dx \geq\int_M(nK-\lambda)\frac{n-1}{n}\lambda v^2\,dx.
\end{equation}

We need show that $\frac{\partial}{\partial \nu}(|\nabla v|^2)\leq
0$ on $\partial M$. Take any $x_0\in \partial M$. If $\nabla
v(x_0)=0$, then it is done. Assume now that $\nabla v(x_0)\not=
0$. Choose a local orthonormal frame $\{e_1, e_2,\cdots, e_{n}\}$
of $M$ about $x_0$ so that $e_n$ is the unit outward normal vector
field near $x_0\in \partial M$ and $\{e_1, e_2,\cdots,
e_{n-1}\}|_{\partial M}$ is a local frame of $\partial M$ about
$x_0$. The existence of such local frame can be justified as the
following. Let $e_n$ be the local unit outward normal vector field
of $\partial M$ about $x_0\in \partial M$ and $\{e_1, \cdots,
e_{n-1}\}$ the local orthonormal frame of $\partial M$ about
$x_0$. By parallel translation along the geodesic
$\gamma(t)=\exp_{x_0}te_n$, we may extend $e_1$, $\cdots$,
$e_{n-1}$ to local vector fields of $M$. Then the extended frame
$\{e_1, e_2,\cdots, e_{n}\}$ is what we need. Note that
$\nabla_{e_n}e_i=0$ for $i\leq n-1$. Since $v|_{\partial M}=0$, we
have $v_i(x_0)=0$ for $i\leq n-1$. Using (\ref{v-d-con})-(\ref{v-d-boundary})
in the following arguments, then we have that at $x_0$,
\begin{eqnarray}\label{gra-v-leq-0}
& &{ }\frac{\partial}{\partial \nu}(|\nabla
v|^2)(x_0)=\sum_{i=1}^n 2v_iv_{in}=2v_nv_{nn}\nonumber\\
& &=2v_n(\Delta^M v -\sum_{i=1}^{n-1}v_{ii}) =2v_n(-\lambda
v-\sum_{i=1}^{n-1}v_{ii})\nonumber\\
&
&=-2v_n\sum_{i=1}^{n-1}v_{ii}=-2v_n\sum_{i=1}^{n-1}(e_ie_iv-\nabla^M_{e_i}e_i
v)\nonumber\\
& &=2v_n\sum_{i=1}^{n-1}\nabla^M_{e_i}e_iv
=2v_n\sum_{i=1}^{n-1}\sum_{j=1}^{n}g(\nabla^M_{e_i}e_i
,e_j)v_j\nonumber\\
& &=2v_n^2\sum_{i=1}^{n-1}g(\nabla^M_{e_i}e_i ,e_n)
=-2v_n^2\sum_{i=1}^{n-1}g(\nabla^M_{e_i}e_n ,e_i)\nonumber\\
& &=-2v_n^2\sum_{i=1}^{n-1}h_{ii}=-2v_n^2(x_0)m(x_0)\nonumber\\
& &\leq 0 \quad\textrm{by the non-negativity of }m,
\end{eqnarray}
where $g(,)$ is the Riemann metric of $M$, $(h_{ij})$ is the
second fundamental form of $\partial M$ with respect to the
outward normal $\nu$ and $m$ is the mean curvature of $\partial M$
with respect to $\nu$. Therefore (\ref{nK-bound}) holds.
\end{proof}

\begin{lemma}       \label{pre-es-d}
Let $v$ be, as the above, the normalized eigenfunction for the
first Dirichlet eigenvalue $\lambda$. Then $v$ satisfies the
following
\begin{equation}                        \label{basic3-d}
\frac{\left |\nabla v\right |^2}{b^2-v^2} \leq\lambda,
\end{equation}
where $b>1$ is an arbitrary constant.
\end{lemma}

\begin{proof}\quad
Consider the function
\begin{equation}                \label{p-of-x-def}
P(x)=|\nabla v|^2+Av^2,
\end{equation}
where $A=\lambda (1 +\epsilon)$ for small $\epsilon>0$. Function
$P$ must achieve its maximum at some point $x_0\in M$. We claim
that
\begin{equation}\label{gra-p-eq-0}
\nabla P(x_0)=0.
\end{equation}
If $x_0\in M\backslash \partial M$, (\ref{gra-p-eq-0}) is
obviously true. Suppose that $x_0\in\partial M$. Take the same
local orthonormal frame $\{e_1, e_2,\cdots, e_{n}\}$ of $M$ about
$x_0$ as in the proof of Lemma \ref{nK-lemma}, where $e_n$ is the
unit outward normal vector field near $x_0\in\partial M$, $\{e_1,
e_2,\cdots, e_{n-1}\}|_{\partial M}$ is a local frame of $\partial
M$ about $x_0$ and $\nabla_{e_n}e_i=0$ for $i\leq n-1$. Since
$v|_{\partial M}=0$, we have $v_i(x_0)=0$ for $i\leq n-1$.
$P(x_0)$ is the maximum implies that
\begin{equation}\label{p-i-eq-0}
P_i(x_0)= 0\qquad \textrm{for }i\leq n-1
\end{equation}
and
\begin{equation}\label{p-n-geq-0}
P_n(x_0)\geq 0.
\end{equation}
Using argument in proving (\ref{gra-v-leq-0}) and the non-negativity
of the mean curvature $m$ of $\partial M$ with respect to the outward
normal, we get
\[
\nabla_{e_n}(|\nabla v|^2)(x_0)\leq 0.
\]
Noticing that $v|_{\partial M}=0$, we have
\begin{equation}\label{p-n-leq-0}
P_n(x_0)=\nabla_{e_n}(|\nabla v|^2)(x_0)+2Av(x_0)v_n(x_0) )\leq 0.
\end{equation}
Now (\ref{p-i-eq-0}), (\ref{p-n-geq-0}) and (\ref{p-n-leq-0})
imply that $P_n(x_0)=0$ and $\nabla P(x_0)=0$.

Thus (\ref{gra-p-eq-0}) holds, no matter $x_0\not\in\partial M$ or
$x_0\in\partial M$. By (\ref{gra-p-eq-0}) and the Maximum
Principle, we have
\begin{equation}\label{max-prin}
\nabla P(x_0)=0 \qquad \textrm{and}\qquad \Delta P(x_0)\leq 0.
\end{equation}

We are going to show further that $\nabla v(x_0)=0$. If on the
contrary, $\nabla v (x_0)\not=0$, then we rotate the local
orthonormal frame about $x_0$ such that
\[
|v_1(x_0)|=|\nabla v(x_0)|\not=0\qquad \textrm{and}\qquad
v_i(x_0)=0,\quad i\geq2.
\]
From (\ref{max-prin}) we have at $x_0$,
\[
0=\frac12\nabla_{i} P=\sum_{j=1}^{n}v_jv_{ji}+Avv_i,
\]
\begin{equation}                \label{d1}
v_{11}=-Av \qquad \textrm{and}\qquad v_{1i}=0\quad i\geq 2,
\end{equation}
and
\begin{eqnarray}
&0 &\geq\frac12 \Delta P(x_0) =\sum_{i, j=1}^{n}\left(v_{ji}v_{ji}+v_{j}v_{jii}+Av_{i}v_{i}+Avv_{ii}\right)\nonumber\\
&{} &=\sum_{i, j=1}^{n}\left(v_{ji}^2+v_j(v_{ii})_j +\textrm{R}_{ji}v_{j}v_{i}+Av_{ii}^2 +Av v_{ii}\right)\nonumber\\
&{} &=\sum_{i, j=1}^{n}v_{ji}^2+\nabla v\nabla(\Delta v) +\textrm{Ric}(\nabla v,\nabla v)+A|\nabla v|^2 +Av\Delta v\nonumber\\
&{} &\geq v_{11}^2+\nabla v\nabla(\Delta v) + (n-1)K|\nabla v|^2+A|\nabla v|^2 +Av\Delta v\nonumber\\
&{} &=(-Av)^2-\lambda |\nabla v|^2+ (n-1)K|\nabla v|^2+A|\nabla v|^2 -\lambda Av^2\nonumber\\
&{} &=(A-\lambda + (n-1)K)|\nabla v|^2+Av^2(A-\lambda),\nonumber
\end{eqnarray}
where we have used (\ref{d1}) and (\ref{ricci-bound}). Therefore
at $x_0$,
\begin{equation}
0\geq (A-\lambda)|\nabla v|^2+A(A-\lambda)v^2, \label{d2}
\end{equation}
that is,
\[
|\nabla v(x_0)|^2+ \lambda(1+\epsilon)v(x_0)^2\leq 0.
\]
Thus $\nabla v(x_0) = 0$. This contradicts $\nabla v(x_0)\not= 0$.

Therefore in any case, if $P$ achieves its maximum at a point
$x_0$, then $\nabla v(x_0)=0$. Thus at $x_0$
\[
P(x_0)=|\nabla v(x_0)|^2 +Av(x_0)^2=Av(x_0)^2\leq A.
\]
and at all $x\in M$
\[
|\nabla v(x)|^2+Av(x)^2=P(x)\leq P(x_0)\leq A.
\]
Letting $\epsilon \rightarrow 0$ in the above inequality, the
estimate (\ref{basic3-d}) follows.
\end{proof}

We want to improve the upper bound in (\ref{basic3-d}) further and
proceed in the following way.

Define a function $F$ by
\[
Z(t)=\max_{x\in M,t=\sin^{-1} \left(v\left(x\right)/b\right)}
\frac{\left |\nabla v\right |^2}{b^2-v^2}/\lambda.
\]
The estimate in (\ref{basic3-d}) becomes
\begin{equation}                                \label{basic5-d}
Z(t)\leq 1\qquad \textrm{on } [0, \sin^{-1}(1/b)]
\end{equation}
For convenience, in this paper we let
\begin{equation}\label{alpha-delta}
\alpha =
\frac12(n-1)K\qquad\textrm{and}\qquad\delta=\alpha/\lambda.
\end{equation}
By (\ref{nK-bound}) we have
\begin{equation}\label{delta-bound}
\delta\leq \frac{n-1}{2n}.
\end{equation}
We have the following conditions for the test function $Z$.

\begin{theorem}                                                     \label{thm-barrier-d}
If the function $z:\ [0,\,\sin^{-1} (1/b)]\mapsto \mathbf{R}^1$
satisfies the following
\begin{enumerate}
 \item $z(t)\geq Z(t) \qquad t\in [0,  \sin^{-1}(1/b)]$,
 \item there exists some $x_0\in M$
       such that at point $t_0=\sin^{-1} (v(x_0)/b)$ \linebreak
       $z(t_0)=Z(t_0)$,
 \item $z(t_0)>0$,
 \item $z$ extends to a smooth even function, and
 \item $z'(t_0)\sin t_0\geq 0$,
\end{enumerate}
then we have the following
\begin{equation}                            \label{barrier-eq-d}
0\leq\frac12z''(t_0)\cos^2t_0 -z'(t_0)\cos t_0\sin t_0 -z(t_0)+
1-2\delta \cos^2t_0.
\end{equation}
\end{theorem}

\begin{proof}\quad
Define
\[ J(x)=\left\{ \frac{\left |\nabla v\right |^2}{b^2-v^2}
-\lambda z \right\}\cos^2t,
\]
where $t=\sin^{-1}(v(x)/b)$. Then
\[ J(x)\leq 0\quad\textrm{for } x\in M
\qquad \textrm{and} \qquad J(x_0)=0.
\]
So $J(x_0)$ is the maximum of $J$ on $M$.
If $\nabla v(x_0)=0$, then
\[ 0=J(x_0)=-\lambda z\cos^2 t.
\]
This contradicts the Condition 3 in the theorem. Therefore
\[ \nabla v(x_0)\not=0.
\]

We claim that
\begin{equation}\label{gra-j-eq-0}
\nabla J(x_0)=0.
\end{equation}
If $x_0\in M\backslash \partial M$, (\ref{gra-j-eq-0}) is
obviously true. Suppose that $x_0\in\partial M$. Take the same
local orthonormal frame $\{e_1, e_2,\cdots, e_{n}\}$ of $M$ about
$x_0$ as in the proof of Lemma \ref{nK-lemma}, where $e_n$ is the
unit outward normal vector field near $x_0\in\partial M$, $\{e_1,
e_2,\cdots, e_{n-1}\}|_{\partial M}$ is a local frame of $\partial
M$ about $x_0$ and $\nabla_{e_n}e_i=0$ for $i\leq n-1$. Since
$v|_{\partial M}=0$, we have $v_i(x_0)=0$ for $i\leq n-1$.
$J(x_0)$ is the maximum implies that
\begin{equation}\label{j-i-eq-0}
J_i(x_0)= 0\qquad \textrm{for }i\leq n-1
\end{equation}
and
\begin{equation}\label{j-n-geq-0}
J_n(x_0)\geq 0.
\end{equation}
Using argument in proving (\ref{gra-v-leq-0}) and the non-negativity
of the mean curvature $m$ of $\partial M$ with respect to the outward
normal, we get
\[
\left(|\nabla v|^2\right)_n\Big|_{x_0}\leq 0.
\]
The Dirichlet condition $v(x_0)=0$ implies that
$t(x_0)=0$ and $z'(t(x_0))=z'(0)=0$, since by the Condition 4 in the theorem
$z$ extends to a smooth even function. Therefore
\begin{equation}\label{j-n-leq-0}
J_n(x_0)=\frac{1}{b^2}\left(|\nabla
v|^2\right)_n-\lambda\cos t[z' \cos t -2z\sin t]t_n\Big|_{x_0}
 =\frac{1}{b^2}\left(|\nabla v|^2\right)_n\Big|_{x_0}
 \leq 0.
\end{equation}

Now (\ref{j-i-eq-0}), (\ref{j-n-geq-0}) and (\ref{j-n-leq-0})
imply (\ref{gra-j-eq-0}).

Thus (\ref{gra-j-eq-0}) holds, no matter $x_0\not\in\partial M$ or
$x_0\in\partial M$. By (\ref{gra-j-eq-0}) and the Maximum
Principle, we have
\begin{equation}                                                \label{es1}
\nabla J(x_0)=0\qquad \textrm{and}\qquad \Delta J(x_0)\leq 0.
\end{equation}
$J(x)$ can be rewritten as
\[  J(x)=\frac{1}{b^2}|\nabla v|^2-\lambda z\cos^2t.
\]
Thus (\ref{es1}) is equivalent to
\begin{equation}                                              \label{es2}
\frac{2}{b^2}\sum_{i}v_iv_{ij}\Big|_{x_0}=\lambda\cos t[z' \cos t
-2z\sin t]t_j\Big|_{x_0}
\end{equation}
and
\begin{eqnarray}                                              \label{es3}
0&\geq&\frac{2}{b^2}\sum_{i,j}v_{ij}^2+\frac{2}{b^2}\sum_{i,j}v_iv_{ijj}
 -\lambda (z''|\nabla t|^2+z'\Delta t)\cos^2t \\
 & &+4\lambda z'\cos t\sin t |\nabla t|^2 -
\lambda z\Delta\cos^2t\Big|_{x_0}.\nonumber
\end{eqnarray}
Rotate the frame so that
\[
|v_1(x_0)|=|\nabla v(x_0)|\not=0\qquad \textrm{and}\qquad
v_i(x_0)=0,\quad i\geq2.
\]
Then (\ref{es2}) implies
\begin{equation}                                             \label{es4}
v_{11}\Big|_{x_0}=\frac{\lambda b}{2}(z'\cos t-2z\sin t)
\Big|_{x_0}\quad\text{and}\quad v_{1i} \Big|_{x_0}=0\ \text{for }
i\geq2.
\end{equation}
Now we have
\begin{eqnarray}
|\nabla v|^2
\Big|_{x_0}&=&\lambda b^2z\cos^2t\Big|_{x_0},\nonumber\\
 |\nabla t|^2
\Big|_{x_0}&=&\frac{|\nabla v|^2}{b^2-v^2}=\lambda z
\Big|_{x_0},\nonumber\\
\frac{\Delta v}{b}\Big|_{x_0} &=&\Delta \sin t =\cos t\Delta
t-\sin t |\nabla t|^2
\Big|_{x_0},\nonumber\\
\Delta t\Big|_{x_0}&=&\frac{1}{\cos t}(\sin t|\nabla
t|^2+\frac{\Delta v}{b})
\nonumber\\
 &=&\frac{1}{\cos t} [ \lambda z\sin t-\frac{\lambda}{b}v] \Big|_{x_0}, \quad\textrm{and}
\nonumber\\
\Delta\cos^2t\Big|_{x_0}&=&\Delta \left(1-\frac{v^2}{b^2}\right)
 =-\frac{2}{b^2}|\nabla v|^2-\frac{2}{b^2}v\Delta v
\nonumber\\
  &=&-2\lambda z\cos^2t+\frac{2}{b^2}\lambda v^2\Big|_{x_0}. \nonumber
\end{eqnarray}
Therefore,
\begin{eqnarray}
& {}&
\frac{2}{b^2}\sum_{i,j}v_{ij}^2\Big|_{x_0}\geq\frac{2}{b^2}v_{11}^2
\nonumber\\
& {}& =\frac{\lambda ^2}{2}(z')^2\cos^2t-2\lambda ^2zz'\cos t\sin
t
      +2\lambda ^2z^2\sin^2 t\Big|_{x_0}\nonumber,
\end{eqnarray}
\begin{eqnarray}
\frac{2}{b^2}\sum_{i,j}v_iv_{ijj}\Big|_{x_0}
&=&\frac{2}{b^2}\left(\nabla v\,\nabla
       (\Delta v)+\textrm{Ric}(\nabla v,\nabla v)\right)\nonumber\\
&\geq& \frac{2}{b^2}(\nabla v\,\nabla (\Delta v)+(n-1)K|\nabla v|^2)\nonumber\\
 &=&-2\lambda^2z\cos^2t+4\alpha \lambda z\cos^2t\Big|_{x_0},\nonumber
\end{eqnarray}
\begin{eqnarray}
&{}&  -\lambda (z''|\nabla t|^2+ z'\Delta t)\cos^2t\Big|_{x_0}\nonumber\\
&{}&=-\lambda^2 zz''\cos^2t-
\lambda^2zz'\cos t\sin t\nonumber\\
&{}&{ }+\frac{1}{b}\lambda^2z'v\cos t \Big|_{x_0},\nonumber
\end{eqnarray}
and
\begin{eqnarray}
&{}&4\lambda z'\cos t\sin t|\nabla t|^2-\lambda z\Delta
\cos^2t\Big|_{x_0}
\nonumber\\
&{}&=4\lambda^2zz'\cos t\sin
t+2\lambda^2z^2\cos^2t-\frac{2}{b}\lambda^2zv\sin t
\Big|_{x_0}.\nonumber
\end{eqnarray}
Putting these results into (\ref{es3}) we get
\begin{eqnarray}                                                   \label{es5}
0&\geq&-\lambda^2zz''\cos^2t+ \frac{\lambda^2}{2}(z')^2\cos^2t
+\lambda^2z'\cos t\left(z\sin t +\sin t\right)
  \nonumber\\
 & & {}+2\lambda^2z^2-2\lambda^2z +4\alpha \lambda z\cos^2t
 \Big|_{x_0},
\end{eqnarray}
where we used (\ref{es4}). Now
\begin{equation}                                                    \label{es6}
z(t_0)>0,
\end{equation}
by the Condition 3 in the theorem. Dividing two sides of (\ref{es5})
by $2\lambda^2z\Big|_{x_0}$, we have
\begin{eqnarray}
0&\geq&-\frac12z''(t_0)\cos^2t_0 +\frac12z'(t_0)\cos t_0\left(\sin t_0+\frac{\sin t_0}{z(t_0)}\right) +z(t_0) \nonumber\\
 & & {}  -1 +2\delta \cos^2t_0+\frac{1}{4z(t_0)}(z'(t_0))^2\cos^2t_0.\nonumber
\end{eqnarray}
Therefore,
\begin{eqnarray}
0&\geq&-\frac12z''(t_0)\cos^2t_0 + z'(t_0)\cos t_0\sin t_0+z(t_0) -1 +2\delta \cos^2t_0\nonumber\\
 & & {}+\frac{1}{4z(t_0)}(z'(t_0))^2\cos^2t_0+\frac12z'(t_0)\sin t_0\cos
 t_0[\frac{1}{z(t_0)}-1].                               \label{es6.1}
\end{eqnarray}
Conditions 1, 2 and 5 in the theorem imply that
$0<z(t_0)=Z(t_0)\leq 1$ and $z'(t_0)\sin t_0\geq 0$. Therefore the last
two terms in (\ref{es6.1}) are nonnegative and
(\ref{barrier-eq-d}) follows.
\end{proof}

%%%%%%%%%Proof of Main Teorems%%%%%%%%%%%%%%%%%%%%
\section{Proof of Main Result}\label{sec-proof}

\begin{proof}[Proof of Theorem \ref{main-thm}]\quad
Let
\begin{equation}                                \label{z-def}
z(t)=1+\delta\xi(t),
\end{equation}
where $\xi$ is the functions defined by (\ref{xi-def}) in Lemma
\ref{xi-lemma}. We claim that
\begin{equation}                        \label{4.1-d}
Z(t)\leq z(t)\qquad \textrm{for }t\in  [0, \sin^{-1}(1/b)].
\end{equation}
Lemma \ref{xi-lemma} implies that for $t\in [0, \sin^{-1}(1/b)] $,
we have the following
\begin{eqnarray}
& &{}\frac{1}{2}z''\cos ^2t-z'\cos t\sin t-z
    =-1+ 2\delta\cos^2t,               \label{z-eq}\\
& &{}z'(t)\sin t\geq 0,\qquad \label{z'-geq0}\\
& &{}z\textrm{ is a smooth even function},\\
& &{}0<1-(\frac{\pi^2}{4}-1)\frac{n-1}{2n}\leq
1-(\frac{\pi^2}{4}-1)\delta=z(0)\leq z(t),\quad\textrm{and}\label{z-min}\\
& &{}z(t) \leq z(\frac{\pi}{2})=1. \label{z-max}
\end{eqnarray}

Let $P\in\mathbf{R}^1$ and $t_0\in [0,\sin^{-1}(1/b)]$ such that
\[ P=\max_{t\in [0,\sin^{-1}(1/b)]}\left(Z(t)-z(t)\right)=Z(t_0)-z(t_0).
\]
Thus
\begin{equation}\label{4.2}
Z(t)\leq z(t)+P\quad \textrm{for }t\in
[0,\sin^{-1}(1/b)]\qquad\textrm{and}\qquad Z(t_0)=z(t_0)+P.
\end{equation}
Suppose that $P>0$. Then $z+P$ satisfies the conditions in
Theorem \ref{thm-barrier-d} and therefore satisfies (\ref{barrier-eq-d}). So we have
\begin{eqnarray}
&{}&z(t_0)+P=Z(t_0)\nonumber\\
&{}&\leq  \frac12(z+P)''(t_0)\cos^2 t_0-(z+P)'(t_0)\cos t_0
 \sin t_0+1-2\delta \cos^2 t_0\nonumber\\
&{}&=\frac12z''(t_0)\cos^2t_0-z'(t_0)\cos t_0\sin
t_0+1-2\delta \cos^2 t_0\nonumber\\
&{}&=z(t_0).\nonumber
\end{eqnarray}
This contradicts the assumption $P>0$. Thus $P\leq 0$ and
(\ref{4.1-d}) must hold. That means
\begin{equation}\label{4.3-d}
\sqrt{\lambda}\geq \frac{|\nabla t|}{\sqrt{z(t)}}.
\end{equation}

Take $q_1$ on $M$ such that $v(q_1)=1 =\sup_M v$ and and $q_2\in
\partial M$ such that distance $d(q_1, q_2) = \textrm{ distance }d(q_1, \partial M)$.
Let $L$ be the minimum geodesic segment between $q_1$ and $q_2$.
We integrate both sides of (\ref{4.3-d}) along $L$ and change
variable and let $b\rightarrow 1$. Let
$\tilde{d}$ be the diameter of the largest interior ball in $M$.
Then
\begin{equation}\label{4.4-d}
\sqrt{\lambda}\,\frac{\tilde{d}}{2}\geq\int_{L}\,\frac{|\nabla
t|}{\sqrt{z(t)}}\,dl =\int_0^{\frac{\pi}{2}}
\frac{1}{\sqrt{z(t)}}\,dt \geq \frac{\left(\int_0^{\pi/2}\
\,dt\right)^\frac32}{(\int_0^{\pi/2}\ z(t)\,dt)^{\frac12}} \geq
\left( \frac{(\frac{\pi}{2})^3}{\int_0^{\pi/2}\  z(t)\,dt}
\right)^{\frac12}.
\end{equation}
Square the two sides. Then
\[
\lambda \geq \frac{\pi^3}{2(\tilde{d})^2\int_0^{\pi/2} \
z(t)\,dt}.
\]
Now
\[
\int_0^{\frac{\pi}{2}}\ z(t)\,dt=\int_0^{\frac{\pi}{2}}\ [1+
\delta \xi(t)]\,dt=\frac{\pi}{2}(1-\delta),
\]
by (\ref{xi-int}) in Lemma \ref{xi-lemma}. That is,
\[
\lambda \geq
\frac{\pi^2}{(1-\delta)(\tilde{d})^2}\quad\textrm{and}\quad
\lambda \geq \frac12(n-1)K + \frac{\pi^2}{(\tilde{d})^2}.
\]
\end{proof}

We now present a lemma that is used in the proof of Theorem
\ref{main-thm}.

\begin{lemma}                           \label{xi-lemma}
Let
\begin{equation}                        \label{xi-def}
\xi(t)=\frac{\cos^2t+2t\sin t\cos t +t^2-\frac{\pi^2}{4}}{\cos^2t}
\qquad \textrm{on}\quad [-\frac{\pi}{2},\frac{\pi}{2}\,].
\end{equation}
Then the function $\xi$  satisfies the following
\begin{eqnarray}
& &{}\frac{1}{2}\xi''\cos ^2t-\xi'\cos t\sin t-\xi
    =2\cos^2t\quad \textrm{in }(-\frac{\pi}{2},\frac{\pi}{2}\,),          \label{xi-eq}\\
& &{}\xi'\cos t -2\xi\sin t =4t\cos t                      \label{xi-eq2}\\
& &{}\int_0^{\frac{\pi}{2}}\xi(t)\, dt= -\frac{\pi}{2}          \label{xi-int}\\
& &{}1-\frac{\pi^2}{4}=\xi(0)\leq\xi(t)\leq\xi(\pm
\frac{\pi}{2})=0\quad
\textrm{on }[-\frac{\pi}{2},\frac{\pi}{2}\,],                                \nonumber\\
& &{} \xi' \textrm{ is increasing on }
[-\frac{\pi}{2},\frac{\pi}{2}\,] \textrm{ and }
\xi'(\pm \frac{\pi}{2}) =\pm \frac{2\pi}{3},                     \nonumber\\
& &{}\xi'(t)< 0 \textrm{ on }(-\frac{\pi}{2},0)\textrm{ and \ }
\xi'(t)>0 \textrm{ on }(0,\frac{\pi}{2}\,), \nonumber \\
& &{}\xi''(\pm\frac{\pi}{2})=2, \ \xi''(0)=2(3-\frac{\pi^2}{4})
\textrm{ and \ } \xi''(t)> 0 \textrm{ on }
[-\frac{\pi}{2},\frac{\pi}{2}\,], \nonumber\\
& &{}(\frac{\xi'(t)}{t})'>0 \textrm{ on } (0,\pi/2\,)\textrm{ and
\ } 2(3-\frac{\pi^2}{4})\leq \frac{\xi'(t)}{t}\leq \frac43
\textrm{ on } [-\frac{\pi}{2},\frac{\pi}{2}\,],\nonumber \\
& &{}\xi'''(\frac{\pi}{2})=\frac{8\pi}{15}, \xi'''(t)< 0 \textrm{
on }(-\frac{\pi}{2},0) \textrm{ and \ }  \xi'''(t)>0 \textrm{ on
}(0,\frac{\pi}{2}\,). \nonumber
\end{eqnarray}
\end{lemma}
\begin{proof}\quad
For convenience, let $q(t)= \xi'(t)$, i.e.,
\begin{equation}                                       \label{q-def}
q(t) = \xi'(t) = \frac{2(2t\cos t +t^2\sin t +\cos^2 t \sin t
-\frac{\pi^2}{4}\sin t)}{\cos^3 t}.
\end{equation}
Equation (\ref{xi-eq}) and the values $\xi(\pm \frac{\pi}{2})=0$,
$\xi(0)=1-\frac{\pi^2}{4}$ and $\xi'(\pm \frac{\pi}{2}) =\pm
\frac{2\pi}{3}$ can be verified directly from (\ref{xi-def}) and
(\ref{q-def}) .  The values of $\xi''$ at $0$ and $\pm
\frac{\pi}{2}$ can be computed via (\ref{xi-eq}). By
(\ref{xi-eq2}), $(\xi(t)\cos^2 t)' =4t\cos^2 t$. Therefore
\newline $\xi(t)\cos^2 t=\int_{\frac{\pi}{2}}^t \ 4s\cos^2 s\,ds$,
and
\[
\int_{-\frac{\pi}{2}}^{\frac{\pi}{2}}\
\xi(t)\,dt=2\int_0^{\frac{\pi}{2}}\
\xi(t)\,dt=-8\int_0^{\frac{\pi}{2}}\left( \frac{1}{\cos^2(t)}
\int_t^{\frac{\pi}{2}}\ s\cos^2s\,ds\right)\,dt
\]
\[
=-8\int_0^{\frac{\pi}{2}}\left(\int_0^s\
\frac{1}{\cos^2(t)}\,dt\right)\ s\cos^2s\,ds
=-8\int_0^{\frac{\pi}{2}}\ s\cos s\sin s\,ds=-\pi.
\]
It is easy to see that $q$ and $q'$ satisfy the following
equations
\begin{equation}                                         \label{q-eq}
\frac12 q''\cos t -2q'\sin t -2q\cos t = -4 \sin t,
\end{equation}
and
\begin{equation}                                       \label{q'-eq}
\frac{\cos^2 t}{2(1+\cos^2 t)}(q')''-\frac{2\cos t\sin t}{1+\cos^2
t}(q')'-2(q')=-\frac{4}{1+\cos^2 t}.
\end{equation}
The last equation implies $q'=\xi''$ cannot achieve its
non-positive local minimum at a point in $(-\frac{\pi}{2},
\frac{\pi}{2})$. On the other hand, $\xi''(\pm\frac{\pi}{2})=2$,
by equation (\ref{xi-eq}), $\xi(\pm \frac{\pi}{2})=0$ and
$\xi'(\pm \frac{\pi}{2})=\pm \frac{2\pi}{3}$. Therefore
$\xi''(t)>0$ on $[-\frac{\pi}{2},\frac{\pi}{2}]$ and $\xi'$ is
increasing. Since $\xi'(t)=0$, we have $\xi'(t)< 0$ on
$(-\frac{\pi}{2},0)$ and $\xi'(t)>0$ on $(0,\frac{\pi}{2})$.
%satisfies the equation
%\[
%\frac12 h'' t\cos t + (\cos t - 2t \sin t) h' -2(\sin t + t\cos t)
%h = -4\sin t.
%\]
Similarly, from the equation
\begin{eqnarray}                                           \label{q''-eq}
&\frac{\cos^2 t}{2(1+\cos^2 t)}(q'')'' -\frac{\cos t\sin t
(3+2\cos^2 t)}{(1+\cos^2 t)^2}(q'')' -\frac{2(5\cos^2 t+\cos^4 t)}{(1+\cos^2 t)^2}(q'') \nonumber\\
&=-\frac{8\cos t\sin t}{(1+\cos^2 t)^2}
\end{eqnarray}
we get the results in the last line of the lemma.

Set $h(t)=\xi''(t)t-\xi'(t)$. Then $h(0)=0$  and $h'(t)=
\xi'''(t)t>0$ in $(0,\frac{\pi}{2})$. Therefore
$(\frac{\xi'(t)}{t})'=\frac{h(t)}{t^2}>0$ in $(0,\frac{\pi}{2})$.
Note that $\frac{\xi'(-t)}{-t}= \frac{\xi'(t)}{t}$,
$\frac{\xi'(t)}{t}|_{t=0}=\xi''(0)=2(3-\frac{\pi^2}{4})$ and
$\frac{\xi'(t)}{t}|_{t=\pi/2}=\frac43$. This completes the proof
of the lemma.
\end{proof}

%%%%%%%%%%%%End of Real Paper%%%%%%%%%%%
%
%
%%%%%%%%%%%%%%%%%References%%%%%%%%%%%%%%%%%%%
%\nocite{*}
%\bibliographystyle{plain}
%\bibliography{GapBound} %GapBound.bib, for using BibTeX
%\input{lbound.bbl}

Department of mathematics, Utah Valley State College, Orem, Utah 84058

\textit {E-mail address}: \texttt{lingju@uvsc.edu}
\end{document}